\documentclass[10pt]{article}
\usepackage{mathrsfs}
\usepackage{amsfonts}
\usepackage{amssymb}  
\usepackage{amsmath,amscd,graphicx}
\usepackage{paralist}  

       \newtheorem{lemma}{\bf Lemma}[section]
       \newtheorem{theorem}[lemma]{\bf Theorem}

       \newtheorem{remark}{\bf Remark}[section]

       \numberwithin{equation}{section}

\begin{document}

\title{{\LARGE \textbf{Blowup of Smooth Solutions to the Navier-Stokes Equations for Compressible Isothermal Fluids}}
 \footnotetext{\small E-mail address: dudp954@nenu.edu.cn(D. Du); lijy645@yahoo.com.cn(J. Li); zhangkj201@nenu.edu.cn(K. Zhang)} }

\author{{Dapeng Du, Jingyu Li and  Kaijun Zhang }\\[2mm]
\normalsize\it   School of Mathematics and Statistics, Northeast Normal University, \\
\normalsize\it   Changchun 130024, P.R.China }

 \date{}

\maketitle

\noindent\textbf{Abstract.} {\small It is shown that the one-dimensional or two-dimensional radially symmetric isothermal compressible
Navier-Stokes system has no non-trivial global smooth solutions if the initial density is compactly supported. This result is a generalization of Xin's work \cite{Xin98} to the isothermal case.
}\\





\section{Introduction}
In this paper we consider the blow-up problem for the Navier-Stokes
equations of compressible isothermal fluids
\begin{equation}\label{NS}
\left\{ \begin{aligned}
&\partial_t\rho+\text{div}(\rho u)=0,\\
&\partial_t(\rho u)+\text{div}(\rho u\otimes u)-\mu\Delta u-(\lambda+\mu)\nabla\text{div}u+\nabla P(\rho)=0,
\end{aligned}\right.
\end{equation}
with initial data
\begin{equation}\label{initial}
\rho(x,0)=\rho_0(x),\quad u(x,0)=u_0(x),\quad x\in\mathbb{R}^n,\quad 1\leq n\leq2.
\end{equation}
This system describes the motion of viscid gas. $u\in \mathbb{R}^n$, $\rho$, $\mu$ and $\lambda$ denote the velocity, the density,
the coefficient of viscosity and the second coefficient of viscosity, respectively. The pressure $P(\rho)$ satisfies
\begin{equation}\label{pressure}
P(\rho)=a\rho
\end{equation} for a positive constant $a$.

\par Mathematically, system \eqref{NS} is a quasilinear transport-parabolic system. There are many important questions
regarding this system. One of them is whether the solutions will develop singularities from smooth initial data. To our best knowledge,
Xin \cite{Xin98} first proved that when the initial densities are compactly supported, any smooth solutions to the compressible Navier-Stokes equations for non-barotropic flows in the absence of heat conduction will blow up in finite time for any spatial dimension, and this feature also holds for the isentropic flows (i.e. $P(\rho)=a\rho^\gamma$ with $\gamma>1$) in one dimensional case. Cho \& Jin \cite{CJ06} extended Xin's work \cite{Xin98} to the case of fluids with positive heat conduction. Recently, under an additional assumption that one of the components of initial momentum is not zero, Tan \& Wang \cite{TW10} gave a much simpler proof of the result of Cho \& Jin \cite{CJ06}. If the initial data do not have compact support but rapidly decrease, for $n\geq3$ and $\gamma\geq\frac{2n}{n+2}$, Rozanova \cite{R08} proved that any smooth solutions to the compressible Navier-Stokes equations for the non-barotropic flows with positive heat conduction still blow up in finite time.

\par It is worth mentioning that very recently, Huang, Li, Luo \& Xin \cite{HLLX} proved that any smooth radially symmetric solutions to the two-dimensional isentropic fluids with compactly supported initial densities blow up. And they also proved in \cite{HLX2} an interesting result that, under some compatibility condition for the initial data, the small smooth solutions to the three-dimensional isentropic compressible Navier-Stokes equations with compactly supported initial densities are globally well-posed.

\par Besides the non-barotropic ones and the isentropic ones, the isothermal flows are also of great physical significance. From the view point of mathematical structures, the isothermal case could be looked as the endpoint of the isentropic case. Roughly speaking, we show that if the initial data are nontrivial and radially symmetric, and the initial densities are compactly supported, then the smooth solutions to  \eqref{NS} will blow up in finite time. More precisely, we prove that:
\begin{theorem}\label{thm}
Assume that
\begin{equation}\label{coeff}
\mu>0,\quad \lambda+\frac{2}{n}\mu>0,
\end{equation}
and that $\rho_0(x)$ has compact support. Suppose that one of the following conditions holds true:

\begin{enumerate}[(i)]

\item $n=1$;

\item $n=2$, and the initial data are spherically symmetric, i.e. $\rho=\rho(|x|,t),u=\bar{u}(|x|,t)x/|x|$.

\end{enumerate}
Then the solution $(\rho,u)\in C^1([0,T], H^m(\mathbb{R}^n))$ $(m>2)$ to system \eqref{NS} with nontrivial initial density will blow up in finite time.

\end{theorem}

\begin{remark}
One key assumption in Theorem \ref{thm} is that the initial density is compactly supported. This assumption could be relaxed to fast decay at spatial infinity like the paper \cite{R08}. But it seems pretty hard to get some kind of blow-up result for the initial density with positive lower bound.
\end{remark}

\begin{remark}
Our method also works for the isentropic case and non-barotropic case. The proofs essentially are the same. One reason is that all these systems have the similar energy conservation and the same momentum equation.
\end{remark}

\par To prove Theorem \ref{thm}, we use the so-called function method which is typical in the theory of blowup. The idea is to integrate the equation suitably and then try to deduce blowup. For the current theorem, the proof roughly goes as follows. First we show that the density keeps being compactly supported, which is similar to that of \cite{Xin98}. Then we integrate the momentum equation by the weight $x$, which gives an integral identity. It consists of three terms. One important observation is that, as time grows, two terms keep bounded, one term grows linearly. This fact implies no nontrivial global smooth solution exists. The details will be presented in Section 2.


\section{Proof of Theorem \ref{thm}}\label{proof of thm}

~~~~In this section we give the proof of Theorem \ref{thm}, which consists of four steps.
\par\emph{Step 1.} The density $\rho$ is compactly supported all the time.
\par Because the initial density $\rho_0(x)$ has compact support, there exists a constant $R>0$ such that
\begin{equation}\label{compact}
\text{supp}\rho_0(x)\subseteq B_R.
\end{equation}
Denote by $x(t;\bar{x})$ the particle path starting from $\bar{x}$, i.e. $x(t;\bar{x})$ satisfies the following equation
\begin{equation}\label{particle}
\left\{ \begin{aligned}
&\frac{dx}{dt}=u(x,t),\\
&x(0)=\bar{x},
\end{aligned}\right.
\end{equation}
where $u$ is the velocity. Denote by $\Omega_t$ the closed region that is the image of $B_R$ under the flow map \eqref{particle},
\begin{equation*}
\Omega_t:=\left\{(x,t)\big|x=x(t;\bar{x}), \bar{x}\in B_R\right\}.
\end{equation*}
Note that on the particle path $x=x(t;\bar{x})$, the density $\rho$ satisfies a homogeneous ordinary differential equation. Thus, if $\rho$ is zero at some initial position $\bar{x}$, then on the particle path starting from $\bar{x}$, $\rho$  will be zero all the time. By \eqref{compact}, we see that
\begin{equation*}
\rho\equiv0 \quad \text{in }\Omega_t^c.
\end{equation*}
Subsequently, using the momentum equation of \eqref{NS}, we get
\begin{equation}\label{velocity}
\mu\Delta u+(\lambda+\mu)\nabla\text{div}u=0 \quad \text{in }\Omega_t^c.
\end{equation}
In the one-dimensional case, owing to \eqref{coeff}, we deduce from \eqref{velocity} that
\begin{equation*}
u_{xx}=0 \quad \text{in }\Omega_t^c.
\end{equation*}
This in combination with the condition $u(\cdot,t)\in H^m(\mathbb{R}^1)$ implies that
\begin{equation*}
u(x,t)\equiv0 \quad \text{in }\Omega_t^c.
\end{equation*}
Therefore, in view of \eqref{particle}, we derive
\begin{equation*}
\Omega_t\equiv\Omega_0\subseteq B_R.
\end{equation*}
In the two-dimensional case, because $u(x,t)=\dfrac{x}{|x|}\bar{u}(|x|,t)$ for some radially symmetric function $\bar{u}$, we have from \eqref{velocity} that
\begin{equation*}
\bar{u}_{rr}+\left(\frac{\bar{u}}{r}\right)_r=0 \quad \text{in }\Omega_t^c.
\end{equation*}
Here $r:=|x|$. Using the condition $u(\cdot,t)\in H^m(\mathbb{R}^2)$ again, it follows that
\begin{equation*}
\bar{u}_r+\frac{\bar{u}}{r}=0 \quad \text{in }\Omega_t^c.
\end{equation*}
One can easily compute that the general solution to this equation is
\begin{equation*}
\bar{u}(r,t)\equiv\frac{C(t)}{r},
\end{equation*}
where $C(t)$ is a constant that only depends on $t$. Because $u\in C([0,T],L^2(\mathbb{R}^2))$, we get for a.e. $t\in[0,T]$,
\begin{equation*}
\int_{\omega_2}\int_0^{+\infty}|\bar{u}(r,t)|^2rdrdS<+\infty,
\end{equation*}
which gives
\begin{equation}\label{comu}
\bar{u}(r,t)\equiv0 \text{  in  } \Omega_t^c.
\end{equation}
Therefore, $\Omega_t\equiv\Omega_0\subseteq B_R$.

\par\emph{Step 2.} Averaging the equation by a weight.
 \par We multiply the momentum equation of \eqref{NS} by weight $x$ and integrate the resulting equation on the whole space to derive
\begin{equation*}
\int_{\mathbb{R}^n}\rho u\cdot xdx-\int_{\mathbb{R}^n}\rho_0 u_0\cdot xdx=\int_0^t\int_{\mathbb{R}^n}\left(\rho|u|^2+na\rho\right)dxd\tau,
\end{equation*}
where we have used the isothermal condition \eqref{pressure}. In view of the fact that
\begin{equation}\label{density}
\int_{\mathbb{R}^n}\rho(x,t)dx=\int_{\mathbb{R}^n}\rho_0(x)dx:=m_0,
\end{equation}
and that $\text{supp}_x\rho(x,t)\subseteq B_R$, we then get
\begin{equation}\label{mom}
\int_{B_R}\rho u\cdot xdx-\int_{B_R}\rho_0 u_0\cdot xdx=\int_0^t\int_{\mathbb{R}^n}\rho|u|^2dxd\tau+nam_0t.
\end{equation}
Obviously, the right-hand side of \eqref{mom} grows linearly. Next we show the left-hand side of \eqref{mom} is bounded. This will give the contradiction with the induction hypothesis that the Navier-Stokes system has a global solution. To get the boundness of the left-hand side of \eqref{mom}, we need to prove energy conservation.

\par\emph{Step 3.} Energy conservation.

\par Multiplying the density equation of \eqref{NS} by $\ln(\rho+\epsilon)$, where $\epsilon>0$ is a small constant, we get
\begin{equation*}
\int_{\mathbb{R}^n}\rho_t\ln(\rho+\epsilon)dx=-\int_{\mathbb{R}^n}\text{div}(\rho u)\ln(\rho+\epsilon)dx.
\end{equation*}
It is easy to see that
\begin{equation*}
\begin{split}
\rho_t\ln(\rho+\epsilon)=[\rho\ln(\rho+\epsilon)]_t-\rho[\ln(\rho+\epsilon)]_t
=[\rho\ln(\rho+\epsilon)]_t-\frac{\rho\rho_t}{\rho+\epsilon}=[\rho\ln(\rho+\epsilon)]_t-\rho_t+\epsilon[\ln(\rho+\epsilon)]_t,
\end{split}\end{equation*}
and that
\begin{equation*}
\begin{split}
\int_{\mathbb{R}^n}\text{div}(\rho u)\ln(\rho+\epsilon)dx&=-\int_{\mathbb{R}^n}\frac{\rho u\cdot\nabla\rho}{\rho+\epsilon}\\
&=-\int_{\mathbb{R}^n}u\cdot\nabla\rho dx+\int_{\mathbb{R}^n}u\cdot\nabla[\epsilon\ln(\rho+\epsilon)] dx\\
&=-\int_{B_R}u\cdot\nabla\rho dx-\int_{B_R}\text{div}u\cdot\epsilon\ln(\rho+\epsilon) dx.
\end{split}
\end{equation*}
Owing to \eqref{density}, we then have
\begin{equation*}
\frac{d}{dt}\int_{B_R}[\rho\ln(\rho+\epsilon)+\epsilon\ln(\rho+\epsilon)]dx
=\int_{B_R}u\cdot\nabla\rho dx+\int_{B_R}\text{div}u\cdot\epsilon\ln(\rho+\epsilon) dx.
\end{equation*}
Let $\epsilon\rightarrow0$, note that $\rho\in C([0,T];H^m(\mathbb{R}^n))$ for $1\leq n\leq2$, the function $\rho\rightarrow\rho\ln\rho$ is continuous with respect to $\rho\in[0,+\infty)$ and that $|\epsilon\ln(\rho+\epsilon)|\leq C|\epsilon\ln\epsilon|\rightarrow0$, it is easy to calculate that
\begin{equation}\label{est-press}
\int_{B_R}\rho\ln\rho dx
=\int_0^t\int_{B_R}u\cdot\nabla\rho dxd\tau+\int_{B_R}\rho_0\ln\rho_0 dx.
\end{equation}
Subsequently, multiplying the momentum equation of \eqref{NS} by $u$ and using the equation of density, we obtain
\begin{equation*}
\frac{d}{2dt}\int_{\mathbb{R}^n}\rho|u|^2dx+\mu\int_{\mathbb{R}^n}|\nabla u|^2dx+(\lambda+\mu)\int_{\mathbb{R}^n}|\text{div} u|^2dx+a\int_{\mathbb{R}^n}u\cdot\nabla\rho dx=0.
\end{equation*}
Subtracting this identity from \eqref{est-press} gives that
\begin{equation*}
\begin{split}
\frac{1}{2}\int_{B_R}\rho|u|^2dx&+a\int_{B_R}\rho\ln\rho dx+\mu\int_0^t\int_{B_R}|\nabla u|^2dxd\tau+(\lambda+\mu)\int_0^t\int_{B_R}|\text{div} u|^2dxd\tau\\&=\frac{1}{2}\int_{B_R}\rho_0|u_0|^2dx+a\int_{B_R}\rho_0\ln\rho_0 dx.
\end{split}
\end{equation*}
Owing to the right-continuity of $\rho\ln\rho$ at $\rho=0$, we get
\begin{equation}\label{gradient}
\int_0^T\int_{B_R}|\nabla u|^2dxd\tau\leq C.
\end{equation}
In the one-dimensional case, by Newton-Leibniz' formula, we get
\begin{equation*}
\|u\cdot x\|^2_{L_x^\infty}\leq R^2\|u\|^2_{L_x^\infty}\leq R^3\int_{-R}^R|u_x|^2dx.
\end{equation*}
In the two-dimensional case, because $u(x,t)=\dfrac{x}{r}\bar{u}(r,t)$, it is easy to calculate that $|\nabla_xu|^2=|\bar{u}|^2/r^2+|\partial_r\bar{u}|^2$. Then
using the Newton-Leibniz' formula and \eqref{comu},
\begin{equation*}
\begin{split}
\|u(x,t)\cdot x\|^2_{L_x^\infty}=\|r\bar{u}(r,t)\|^2_{L_r^\infty}&\leq2R\Big(\int_0^R|\bar{u}(r,t)|^2dr
+\int_0^Rr^2|\partial_r\bar{u}(r,t)|^2dr\Big)\\& \leq2R^2\Big(\int_0^R\frac{|\bar{u}(r,t)|^2}{r^2}\cdot rdr
+\int_0^R|\partial_r\bar{u}(r,t)|^2rdr\Big)\\
&\leq \frac{R^2}{2\pi}\int_{B_R}|\nabla_x u(x,t)|^2dx.
\end{split}
\end{equation*}
Consequently,
\begin{equation}\label{bound}
\Big(\int_{B_R}|\rho u\cdot x|dx\Big)^2\leq m_0^2\|u\cdot x\|^2_{L_x^\infty}\leq C\int_{B_R}|\nabla u(x,t)|^2dx.
\end{equation}

\par\emph{Step 4.} Derivation of the blow-up.

\par By \eqref{mom}, \eqref{bound} and \eqref{gradient}, we have
\begin{equation*}
\begin{split}
a^2m_0^2T^3&\leq2\int_0^T\Big(\int_{B_R}\rho u\cdot xdx\Big)^2dt+2\int_0^T\Big(\int_{B_R}\rho_0 u_0\cdot xdx\Big)^2dt\\&
\leq C+CT.
\end{split}
\end{equation*}
Therefore, it follows that the smooth solutions to \eqref{NS} will blow up in finite time.
\hspace{\stretch{1}}$\square$

\section*{Acknowledgements} Du's work is partially supported by the Chinese NSF (No. 11001043) and the Chinese Postdoctoral Science Foundation (No. 20090460074);
Li's work is partially supported by the Fundamental Research Funds for the Central Universities (No. 10QNJJ001);
Zhang's work is partially supported by the Chinese NSF (No. 11071034) and the Fundamental Research Funds for the Central Universities (No. 111065201).

\bibliographystyle{amsplain}

\end{document}